\documentstyle[amssymb,10pt]{amsart}
\newtheorem{thm}{Theorem}[section]

\theoremstyle{remark}
\newtheorem{remark}{Remark}

\begin{document}
\newcommand{\nc}{\newcommand} \nc{\on}{\operatorname}
\nc{\pa}{\partial}
\nc{\cA}{{\cal A}}\nc{\cB}{{\cal B}}\nc{\cC}{{\cal C}}
\nc{\cE}{{\cal E}}\nc{\cG}{{\cal G}}\nc{\cH}{{\cal H}}
\nc{\cX}{{\cal X}}\nc{\cR}{{\cal R}}\nc{\cL}{{\cal L}}
\nc{\cK}{\mathcal {K}}\nc{\cO}{\mathcal {O}}
\renewcommand{\i}{\on{i}}
\nc{\sh}{\on{sh}}\nc{\Id}{\on{Id}}\nc{\Diff}{\on{Diff}}
\nc{\ad}{\on{ad}}\nc{\Der}{\on{Der}}\nc{\End}{\on{End}}\nc{\res}{\on{res}}
\nc{\Imm}{\on{Im}}\nc{\limm}{\on{lim}}\nc{\Ad}{\on{Ad}}
\nc{\Hol}{\on{Hol}}\nc{\Det}{\on{Det}}
\nc{\de}{\delta}\nc{\si}{\sigma}\nc{\ve}{\varepsilon}
\nc{\al}{\alpha}
\nc{\CC}{{\mathbb C}}\nc{\ZZ}{{\mathbb Z}}\nc{\NN}{{\mathbb N}}
\nc{\AAA}{{\mathbb A}}
\nc{\cF}{{\cal F}}
\nc{\la}{{\lambda}}\nc{\G}{{\mathfrak g}}\nc{\A}{{\mathfrak a}}
\nc{\HH}{{\mathfrak h}}
\nc{\N}{{\mathfrak n}}\nc{\B}{{\mathfrak b}}
\nc{\La}{\Lambda}
\nc{\g}{\gamma}\nc{\eps}{\epsilon}\nc{\wt}{\widetilde}
\nc{\wh}{\widehat}
\nc{\bn}{\begin{equation}}\nc{\en}{\end{equation}}
\nc{\SL}{{\mathfrak{sl}}}
\nc{\vj}{\mathbf{v}}
\nc{\I}{\mathbf{I}}
\nc{\ny}{\nonumber}
\def\r#1{(\ref{#1})}
\let\ds=\displaystyle
\def\sk#1{\left( #1 \right)}
\def\skk#1{\left[ #1 \right]}
\def\br#1{\left\langle #1\right\rangle}

%
%
%

\newcommand{\ldar}[1]{\begin{picture}(10,50)(-5,-25)
\put(0,25){\vector(0,-1){50}}
\put(5,0){\mbox{$#1$}}
\end{picture}}

\newcommand{\lrar}[1]{\begin{picture}(50,10)(-25,-5)
\put(-25,0){\vector(1,0){50}}
\put(0,5){\makebox(0,0)[b]{\mbox{$#1$}}}
\end{picture}}

\newcommand{\luar}[1]{\begin{picture}(10,50)(-5,-25)
\put(0,-25){\vector(0,1){50}}
\put(5,0){\mbox{$#1$}}
\end{picture}}

\title{Basic representations of quantum current algebras in higher genus}

\dedicatory{To the memory of Joseph Donin}

\author{B. Enriquez}\address{IRMA, Universit\'e Louis Pasteur,
Strasbourg, 7, rue Ren\'e Descartes, F-67084 Strasbourg, France}
\author{S. Pakuliak}\address{Laboratory of Theoretical Physics, JINR,
141980 Dubna, Moscow reg., Russia}
\author{V. Rubtsov}\address{D\'epartment de Math\'ematiques,
Universit\'e d'Angers,
2 Bd. Lavoisier, 49045 Angers, France and ITEP Theory Department,
25, Bol. Tcheremushkinskaya, 117259, Moscow, Russia}

\date{\today}
\maketitle

\begin{abstract}
We construct level 1 basic representations of the quantized
current algebras associated to higher genus algebraic curves
using one free field.
We also clarify the relation between the elliptic current
algebras of the papers \cite{EF98} and \cite{JKOS}.
\end{abstract}

\section{Introduction}
The recent decade has been a period of interest to different
quantizations of current algebras associated with algebraic
curves. Special attention was paid to elliptic quantum
current algebras. The first paper where their main ingredients
were introduced is probably
\cite{L95}. The screening currents of a scaling elliptic algebra
appear there for the first time. The algebra of screening currents was later
studied in \cite{KLP-AMS,CKP}.

The first relation between elliptic algebras and quasi-Hopf algebras
(\cite{D90})
appeared in the papers \cite{HGQG,Fr97}. In the first paper, elliptic algebras
are expressed using total currents, in the framework of Drinfeld's
comultiplication. It became later clear that applications to integrable models
require the introduction of half-currents and Gauss coordinates of
$L$-operators. The first paper which used the formalism of $L$-operators
in the context of elliptic current algebras was \cite{FIJKMY}, but
the quasi-Hopf nature of these algebras was not observed there.

The paper \cite{EF98} used the quasi-Hopf constructions of
\cite{HGQG,ER99} and introduced a decomposition of total currents into
half-currents, which were identified with the Gauss coordinates
of $L$-operators (in analogy with \cite{DF}), and the resulting current
algebra was written in terms of $L$-operators and identified with
the central extension of Felder's elliptic quantum group
$E_{\tau,\eta}(\widehat{\mathfrak{sl}}_2)$.

Precise relations between some elliptic current algebras in the quasi-Hopf
framework (starting directly from the $L$-operator formulation) were obtained
in the papers \cite{ABRR97,JKOS97}. The corresponding algebras also
satisfy the dynamical Yang-Baxter equations; in this case, the dynamical shift
not only affects the dynamical variable by an element of the Cartan algebra
like in the usual $L$-operator formalism, but it also changes the elliptic
modulus by a central element.  Such deformations, leading to shifts in
the moduli of curves, are known only in the case of genus one (elliptic curves).

This paper is a new contribution to the study of quantum current algebras
associated with algebraic curves. Namely, we construct level 1 basic
representations for these algebras using one free field. This is
the analogue of
the Frenkel-Jing construction \cite{FJ} (see Frenkel-Kac \cite{FK}
in the classical case). Recall that zero-level
representations of these algebras were studied in \cite{HGQG}. The next step
of our program is the construction of a bundle of level 1 vertex
operators, a natural
connection on this bundle and an analogue of the qKZB discrete
connection. We hope to return to these questions elsewhere.

The paper is organized as follows. In Section 2, we recall the
general construction of the papers \cite{HGQG,ER99} of quantum
current algebras on algebraic curves. We recall the main ingredients
of this construction: the Green function and the function $q(z,w)$,
and we give some examples is low genus. The next section is devoted
to the main result of the paper: the construction of level one basic
representations of the quantum current algebra associated with
$\mathfrak{sl}_2$
and arbitrary algebraic curves. This construction resembles the approach
developed in the papers \cite{Kh,KLP96}. The last section contains a
comparison of the total currents decomposition into half-currents in the
elliptic quantum group $E_{\tau,\eta}(\widehat{\mathfrak{sl}}_2)$ and
in the algebra $U_{q,p}(\widehat{\mathfrak{sl}}_2)$ of the paper \cite{JKOS}.
(As Jimbo explained to us, the decomposition in \cite{EF98} was a starting
point for the paper \cite{JKOS}, in which the algebras have the advantage of
being graded.)
We present some arguments that these algebras are elliptic generalizations
of two presentations of the centrally extended Yangian double given
in the paper \cite{KLP-AMS}. In the conclusion, we discuss
several open problems.

\section{Green kernels and the function $q(z,w)$. Description and examples.}

\subsection{Notations and conventions} Let $\Sigma$ be a smooth connected
complex compact curve and $\omega$ be a nonzero meromorphic differential on
$\Sigma$. Let $S_{\omega} \subset \Sigma$ denote the set of
poles and zeros of $\omega$. If ${\cK}_s$ denotes the local
field at the point $s\in S_{\omega}$ and ${\cK} :=
\oplus_{s\in S_{\omega}}{\cK}_s$, then the ring $R$ of functions regular
outside $S_{\omega}$ can be viewed as a subring in $\cK$. We also choose
a subspace $\Lambda$ in $\cK$ which is a {\it Lagrangian complement}
to $R$ with respect to the bilinear form
$\langle f,g\rangle_{\cK}= \sum_{s\in S_{\omega}}
\on{res}_{s}(fg\omega)$,
commensurable to $\oplus_{s\in S_\omega}\cO_s$ ($\cO_s$ is the local ring at
$s$). Finally, we define a derivation $\partial$ of $\cK$ by
$\partial f = {{df}\over{\omega}}$.

We will adopt here the following notations from our earlier works: $\hbar$ is
a formal variable, $q = e^\hbar$.
The operator $e^{\hbar\partial}$ is an automorphism of ${\cK}[[\hbar]]$
preserving $R[[\hbar]]$.
We will choose a local coordinate $z_s$ on $\Sigma$ at
$s\in S_{\omega}$ and we will set
$V((z)) = \prod_{s\in S_{\omega}}V[[z_s]][z_{s}^{-1}]$
for any vector space $V$; this is a
completion in $V\otimes {\cK}$. The ring $\mathbb C[[z,w]][z^{-1},w^{-1}] =
\prod_{(s,t)\in S_{\omega}\times S_{\omega}}{\mathbb C}
[[z_s,w_t]][z_{s}^{-1},w_{t}^{-1}]$ is a completion
of ${\cK}\otimes{\cK}.$

We will denote by $f\to f^{(21)}$ the permutation of factors in
$\mathbb C[[z,w]][z^{-1},w^{-1}]$. If
$f=\sum_k f'_{k}\otimes f''_{k} \in R\otimes R$, we
set $f(z,w)= \sum_k f'_{k}(z)f''_{k}(w)$, a
complex function in
$(\Sigma \backslash S_{\omega})\times (\Sigma \backslash S_{\omega})$.
In this case $f^{(21)}(z,w) = f(w,z).$
For any $f\in {\cK}$ we will always denote by $(f)_{R}$ (resp., $(f)_{\Lambda}$) the projections of
$f$ on $R$ parallel to $\Lambda$ (resp., on $\Lambda$ parallel to $R$).
Finally, we will choose $r^i, i\geq 0$ and $\lambda_i, i\geq 0$
two dual bases of $R$ and $\Lambda$.
The following notation for the group commutator will be used below:
$$
(a,b):=aba^{-1}b^{-1}
$$
\subsection{Green kernels} Let us define the {\it Green kernel} by
$G(z,w)=\sum_{i\geq 0} r^i (z)\lambda_i (w) \in \mathbb C((z))((w))$.
We summarize its main properties:
\begin{thm} \label{thm21} (see \cite{ER99})
\begin{equation}
((\partial\otimes \on{id})G)(z,w) = -G(z,w)^2 + \gamma,
\end{equation}
where $\gamma \in R^{\otimes 2}$.

Let $\phi,\psi$ be the elements
of $\hbar R^{\otimes 2}[[\hbar]]$ satisfying the differential equations
\begin{equation}
\partial_{\hbar}\psi = (\partial\otimes \on{id})\psi - 1 - \gamma\psi^2, \
\partial_{\hbar}\phi = (\partial\otimes \on{id})\phi - \gamma\psi.
\end{equation}
Then $\psi = -\hbar + o(\hbar)$, $\phi = (1/2) \hbar^2 \gamma + o(\hbar^2)$.
Then the following identity takes place:
\begin{equation}
\sum_{i\in \mathbb N}{{e^{\hbar\partial}-1}\over{\partial}}\ \lambda_i (z)r^i (w) =
-(\phi(\hbar) + \on{log}(1-G^{(21)}\psi(\hbar)))(z,w).
\end{equation}
\end{thm}

The element
$\sum_{i\geq 0}r^i \otimes({{e^{\hbar\partial}-
e^{-\hbar\partial}}\over{\partial}}\lambda_i)_R$
belongs to $S^2(R)[[\hbar]]$. Let us fix an element $\tau$ in
$R^{\otimes 2}[[\hbar]]$ such that
\begin{equation}
\hbar\sk{\tau + \tau^{(21)}} + \sum_{i\geq 0}r^i
\otimes\sk{{{e^{\hbar\partial}- e^{-\hbar\partial}}
\over{\partial}}\lambda_i}_R = 0.
\end{equation}

We define $\delta(z,w) = \sum_i \epsilon^i (z)\epsilon_i (w)$,
where $(\epsilon^i)$
and $(\epsilon_i)$ are dual bases for $\cK$ for $\langle\ ,\ \rangle_{\cK}.$
Then $\delta(z,w) = (G + G^{(21)})(z,w)$. We have $(\pa\otimes \on{id} +
\on{id}\otimes \pa)(\delta(z,w)) = 0$. One can view
the Green kernel as the collection
of expansions in $w$ at the vicinity of $s\in S_{\omega}$ of a rational function on
$\Sigma\times \Sigma$ antisymmetric in $z$ and $w$, regular except the poles where $z(w)=s$ for some
$s\in S_{\omega}$ and with a simple pole on diagonal.

\subsection{The function $q(z,w)$} Let us define a function $q(z,w)\in
{\mathbb C}((w))((z))[[\hbar]]$
 by the formula:
\begin{equation}
q(z,w) = \exp\left(\sum_i \left({{e^{\hbar\partial}-
e^{-\hbar\partial}}\over{\partial}}
\lambda_i\right)
\otimes r^i\right)\ \exp(\hbar\tau)(z,w).
\end{equation}

Let us denote by $\CC[(\Sigma \setminus S_\omega)^2 \setminus \on{diag}]$
the ring of regular functions on $(\Sigma \setminus S_\omega)^2 \setminus
\on{diagonal}$. Theorem \ref{thm21} implies that $q(z,w)$ is the expansion,
for $z$ close to $S_\omega$, of an element $\widetilde q(z,w)\in
\CC[(\Sigma \setminus S_\omega)^2 \setminus \on{diag}][[\hbar]]$, of
the form $1+O(\hbar)$ and such that $\widetilde q(z,w)\widetilde q(w,z)=1$.

We will also use another property of $q(z,w)$. Assume that
$\on{card}(S_\omega) = 1$,
then there exists a series $i(z,w)\in \CC[[z,w]][z^{-1},w^{-1}][[\hbar]]$
of the form $1 + O(\hbar)$, such that
$q(z,w) = i(z,w)(z-q^{\partial}w)/(q^{\partial}z-w)$; we have $i(z,w)i(w,z)=1$.

\subsection{Examples of the function $q$}
\subsubsection{$g=0$, $\omega = dz$} Assume
$\Sigma = {\mathbb C}P^1, \omega = dz, S_{\omega}=\{\infty\}.$ The local
coordinate is $z_{\infty}=z^{-1}$.
We have $R={\mathbb C[[z]]}$ and take as a complement
$\Lambda = z^{-1}{\mathbb C}[[z^{-1}]].$
Then the Green kernel $G(z,w)$ is given by the expansion of
${{1}\over{z-w}}$ for $z\ll w$ and the function
$q(z,w,\hbar) = {{z-w+\hbar}\over{z-w-\hbar}}$ for
$w\ll z.$

\begin{remark}
This example is a particular ($N=1$) case of the
following more general construction (see \cite{ER-MS}):
assume $\Sigma =
{\mathbb C}P^1, \omega = z^{N-1}dz, S_{\omega}=\{\infty\}.$

If $N$ is odd and we take $N=2n+1$ then $R = z^{-n-1}
{\mathbb C}[z^{-1}]$ and $\Lambda = z^{-n}{\mathbb C}[[z]].$
The dual bases are $r^i = z^{-n-i-1}$ and
$\lambda_i = z^{i-n}$ for $i\geq 0$. Then the Green kernel
$G = \sum_{i\geq 0} z^{-n-i-1}\otimes z^{i-n} = {{(zw)^{-n}}\over{z-w}}$ expands at $w=0.$
\end{remark}

\subsubsection{\bf $g=0$, $\omega = {{dz}\over{z}}$} Assume $\Sigma =
{\mathbb C}P^1, \omega = {{dz}\over{z}},
S_{\omega}=\{0,\infty\}.$ We pose $\Lambda = \{(\lambda_0; \lambda_{\infty})\in
{\mathbb C}[[z]]\times
{\mathbb C}[[z^{-1}]]|\lambda_0 (0) +
\lambda_{\infty}(\infty)=0\}.$ Dual bases for $R$ and for $\Lambda$
are $\{r^i = z^i\}, i\in {\mathbb Z}$ and
$\lambda_i = (z^{-i}; 0)$ for $i<0$, $\lambda_i = (0;-z^{-i})$ for
$i>0$ and $\lambda_0 = (1;-1)$. We can compute easily
$$
\exp\left(\sum_i\left({{e^{\hbar\partial}-e^{-\hbar\partial}}\over{\partial}}
\right)r^i\otimes\lambda_i\right)(z,w)
= \left({{e^{\hbar}z-w}\over{z-e^{\hbar}w}};{{e^{\hbar}z-w}\over{z-e^{\hbar}w}}\right).
$$
\subsubsection{$g=1$, $\omega = dz$}\label{2.4.3}
Assume $\Sigma = {\mathbb C}/{\Gamma},$where $\Gamma = \mathbb Z
+\tau\mathbb Z$ for
a fixed complex number $\tau, \Im\tau>0.$ We choose a coordinate $z$ on
$\Sigma$ and $\omega = dz$. In our
case $\cK$ is the completed local field ${\mathbb C}((z))$ at $z=0$ and
$R={\mathbb C}[[z]]$. It is clear that
$R$ is a maximal subring in $\cK$ isotropic with respect to
$\langle f,g\rangle_{\cK}= \on{res}_{z=0}(fg dz).$

Let $\theta$ be the theta-function, defined to be
holomorphic in $\mathbb C$ with only zeros at the points
of the lattice $\Gamma$, such that
\footnote{We denote $\i = \sqrt{-1}$, leaving $i$ for indices.}
$\theta(z+1) = -\theta(z)$, $\theta(z+\tau) = -e^{\i\pi(2z+\tau)}\theta(z)$
and $\theta'(0)=1.$
Define a complementary subspace $\Lambda(\lambda), \lambda \in {\mathbb C}$
to $R$:
\begin{equation}
\Lambda(\lambda) = \oplus_{k\geq 0}{\mathbb C}\cdot\partial^{k}\zeta(z)
e^{2\i\pi\Re(\lambda)z}, \ \lambda \in \Gamma,
\end{equation}
\begin{equation}
\Lambda(\lambda) = \oplus_{k\geq 0}{\mathbb C}\cdot\partial^{k}{{\theta(\lambda + z)}\over{\theta(z)}}, \ \lambda \not\in \Gamma,
\end{equation}
where $\zeta(z) = {\frac{d}{dz}}(\on{log}(\theta(z)))
={{\theta'}\over{\theta}}(z).$
\begin{remark}
 $\Lambda(0)$ is a maximal
isotropic subspace in $\cK$, consisting of all $\Gamma$-periodic
functions, holomorphic on $\mathbb C\backslash \Gamma$ and such that the $\oint_{\gamma_\epsilon}f(z)dz=0$, where
$\epsilon>0$ and the closed cycle $\gamma_{\epsilon}=(\i\epsilon,1+\i\epsilon).$
\end{remark}
The Green kernel $G(z,w)$ in this case is $\zeta(z-w)-\zeta(z)+\zeta(w)$ and $q(z,w,\hbar) =
{{\theta(z-w+\hbar)}\over{\theta(z-w-\hbar)}}.$

\section{Quantum current algebra associated with an algebraic curve}

We define $U_\hbar(\G)$ as the topologically free $\CC[[\hbar]]$-algebra
with generators $K$, $h^+[r]$,
$h^-[\lambda]$, $x^\pm[\eps]$, where $r\in R$, $\lambda\in \Lambda$,
$\eps\in\cK$ and relations: each map $\varphi\mapsto a[\varphi]$
is linear, $K$ is central, and if we set
$h^+(z) = \sum_i h^+[r_i] \lambda_i(z)$,
$h^-(z) = \sum_i h^-[\lambda_i]r_i(z)$,
$$
K^+(z) = e^{\hbar((T+U)h^+)(z)}, \; K^-(z) = e^{\hbar h^-(z)}, \;
x^\pm(z) = \sum_i x[\eps^i]\eps_i(z),  \; (x = e,f),
$$
where $U : \Lambda \to R$ is defined by $U(\lambda) = \langle \tau,
\on{id}\otimes \lambda\rangle$, then relations are
\begin{equation}\label{e-e}
\Big( e(z)e(w) - q(z,w)e(w)e(z) \Big) \alpha(z,w) = 0,
\end{equation}
\begin{equation}\label{f-f}
\Big( f(z)f(w) - q(w,z)f(w)f(z) \Big) \alpha(w,z) = 0,
\end{equation}
\begin{equation}\label{K-e}
K^+(z)e(w)K^+(z)^{-1} = q(z,w) e(w),
\quad  K^-(z)e(w)K^-(z)^{-1} = q(w,e^{-\hbar K\pa}z) e(w),
\end{equation}
\begin{equation}\label{K-f}
K^+(z)f(w)K^+(z)^{-1}=q(z,w)^{-1}f(w),\quad
K^-(z)f(w)K^-(z)^{-1}=q(w,z)^{-1}f(w)
\end{equation}
\begin{equation}\label{Kp-Kp}
(K^+(z),K^+(w))=1
\end{equation}
\begin{equation}\label{Kp-Km}
(K^+(z),K^-(w))=q(z,w)q(z,e^{-\hbar K\pa}w)^{-1},
\end{equation}
\begin{equation}\label{Km-Km}
(K^-(z),K^-(w))=
q(e^{-\hbar K\pa}z,e^{-\hbar K\pa}w) q(z,w)^{-1},
\end{equation}
\begin{equation}\label{e-f}
[e(z),f(w)]={1\over\hbar}
\Big(\delta(z,w)K^+(z)-\delta(z,e^{-\hbar K\pa}w)
K^{-}(w)^{-1} \Big).
\end{equation}
Here $q(z,w)$ is as above, and $\alpha(z,w)$ runs over all elements of
$\CC[[z,w]][z^{-1},w^{-1}][[\hbar]]$ such that
$\alpha(e^{\hbar\partial}w,w)=0$.

For the function
\begin{equation}\label{qYang}
q(z,w)=\frac{z-w+\hbar}{z-w-\hbar}
\end{equation}
this algebra was considered in \cite{Kh} and it coincides with
a centrally extended Yangian double associated to $\SL_2$.
For the function
\begin{equation}\label{qEF}
q(z,w)=\frac{\theta(z-w+\hbar)}{\theta(z-w-\hbar)}
\end{equation}
this algebra was considered in \cite{EF98} and after proper
factorization of the currents $e(z)$ and $f(z)$ was identified in
this paper with the elliptic quantum group constructed in
\cite{Ell-QG}.

\section{Basic representations of $U_\hbar(\G)$}

In this section, we assume that $\on{card}(S_\omega)=1$.
Let us fix a local coordinate
$z$ on $\Sigma$, so that $\cK = \CC((z))$. We define $\alpha_0 =
d\on{log}z/\omega(z)\in\cK$. We define $\wt\cK
:= \cK\oplus \CC \on{log}z$. We then define
$\partial : \wt\cK \to\cK$ as the extension of $\partial$,
taking $\on{log}z$ to $\alpha_0$. We also define
$(\partial^{-1}) : \cK \to \wt\cK$
as the unique linear map, taking $f\in\cK$ such that $\langle f,1\rangle=0$
to the element $g\in\cK$ such that $\partial(g)=f$ and
$\langle g,\alpha_0\rangle=0$,
and taking $\alpha_0$ to $\on{log}z$.
Then $\partial \circ (\partial^{-1}) = \on{id}_{\cK}$ and
$(\on{id} - (\partial^{-1})\circ \partial)(f) = \langle f,\alpha_0 \rangle 1$
for $f\in\cK$.

Let us define a Heisenberg algebra $H$ as the topologically free
$\CC[[\hbar]]$-algebra with generators
$A[\eps]$, $c$ ($\eps\in\cK$) and relations: $\eps\mapsto A[\eps]$
is linear,
$$
[A[r],A[r']] = 0, \; [A[r],A[\lambda]] = \hbar^{-1}\langle
(1-e^{-\hbar\partial})(r), \lambda\rangle,
$$
$$
[A[\lambda],A[\lambda']] = \hbar^{-1}\br{
(e^{-\hbar\partial} \otimes e^{-\hbar\partial}-1)\sk{\sum_i
(T+U)(\lambda_i) \otimes r^i}, \lambda \otimes \lambda'};
$$
$$
[c,A[r]] = \br{ {{e^{-\hbar\partial}-1}\over{\hbar\partial}}\ (r),
\alpha_0},
\quad [c,A[\lambda]] = \br{
{{e^{-\hbar\partial}-1}\over{\hbar\partial}}\ (T+U)((e^{\hbar\partial}
\lambda)_\Lambda), \alpha_0}
$$
($r,r'\in R$, $\lambda,\lambda'\in\Lambda$);
here $T = (e^{\hbar\partial} - e^{-\hbar\partial})/(\hbar\partial)$.
Note that $(e^{-\hbar\partial}
\otimes e^{-\hbar\partial}-1)(\sum_i
(T+U)(\lambda_i) \otimes r^i) \in \wedge^2(R)[[\hbar]]$.

The subalgebra of $H$ generated by the $A[\eps]$
identifies with the Cartan currents
subalgebra of $U_\hbar(\G)$ generated by the $h^+[r]$, $h^-[\lambda]$,
where $K=1$.

Let us set $A_R(z) = \sum_i A[r^i](T+U)(\lambda_i)(z)$,
$A_\Lambda(z) = \sum_i A[\lambda_i] r^i(z)$.

We then set
$$
a_\Lambda(z) := {{\hbar\partial}\over{1-q^{-\partial}}}
(\partial^{-1})(A_\Lambda(z)),
\quad
a_R(z) := {{\hbar\partial}\over{q^\partial-1}}
(\partial^{-1})(A_R(z)).
$$

We set
$$
a(z) := q^\partial a_R(z) + c + a_\Lambda(z), \quad
b(z) := a_R(z) + c + a_\Lambda(z).
$$

Note that $a(z) - b(z) = \hbar A_R(z)$, $a(z) - b(q^\partial z) =
- \hbar A_\Lambda(q^\partial z)$.

If $f(z) = \bar f \on{log}(z) + \sum_{n\in\ZZ} f_i z^i$, we set
$f_+(z) = \sum_{i\geq 0} f_i z^i$ and $f_-(z) = \bar f \on{log}z +
\sum_{i<0} f_i z^i$. We then split $a(z) = a_+(z) + a_-(z)$,
$b(z) = b_+(z) + b_-(z)$.

We define special functions $v(z)$, $v'(z)$, $u_R(z)$ and $u_\Lambda(z)$
as follows.

We first set
$j(z,w) := {{z-q^\partial w}\over{q^{-\partial}z-w}}$.
We have $j(z,w)\in \CC[[z,w]][z^{-1},w^{-1}]
[[\hbar]]$, $j = 1 + O(\hbar)$.

We then set $\phi(z,w):= \Big( {1\over 2}\log(j^{21})_{++}
+ {1\over 2}\log(j^{21})_{--}\Big)(z,w)$,
then
$\phi(z,w)\in \hbar \CC[[z,w]][z^{-1},w^{-1}][[\hbar]]$.

There is a unique map $R^{\otimes 2} \to \cK$, $a(z)b(w)\mapsto
a_-(z) b_+(z)$, which we denote by $f(z,w)\mapsto f(z,w)_{-+|w=z}$.
We then set $\alpha(z) = \phi(z,z)$ and $\beta(w) = \phi(q^{-\pa}w,w)
+ \Big((q^{-\pa}\otimes q^{-\pa}-1)\log q(z,w)\Big)_{-+|w=z}$.

Recall that $\hbar\cK[[\hbar]] = \hbar R[[\hbar]] \oplus
(T+U)(\hbar\Lambda[[\hbar]])$. If $\rho\in
\hbar\cK[[\hbar]]$, we denote by
$\rho = \rho_R + \rho_\Lambda$ the corresponding decomposition.

We then set
$$
u_R(z):= e^{\alpha_R(z)}, \quad
u_\Lambda(z) := e^{-\beta_\Lambda(w)}.
$$
We define $r_0(z)\in\cK$ by $\omega(z) = r_0(z)dz$ and set
$$
v(z) := (z - q^\pa z)/\hbar
, \quad
v'(z) := r_0(z)e^{-\alpha_\Lambda(z) - \beta_R(z)}.
$$

\begin{thm}
The functional relations defining $U_\hbar(\G)$ are satisfied by the
substitutions
$$
K\mapsto 1,\quad K^+(z) \mapsto u_R(z) e^{\hbar A_R(z)},
\quad K^-(z) \mapsto u_\Lambda(z) e^{\hbar A_\Lambda(z)},
$$
$$
e(z) \mapsto  E(z) :=
v(z) \exp(a_-(z)) \exp(a_+(z))
$$
$$
f(z) \mapsto  F(z) :=
v'(z) \exp(-b_-(z)) \exp(-b_+(z)).
$$
These formulas induce a functor $\{$topological $H$-modules such that
$A[1]$ is diagonalizable with eigen\-va\-lues in $\ZZ\} \to
\{$level
$1$ modules over $U_\hbar(\G)\}$. Here topological means that for
any $v\in V$, $\langle A_R(z), z^n\rangle v$ tends
$\hbar$-adically to $0$ when $n\to+\infty$.
\end{thm}

{\em Proof.} The images of $K^\pm(z)$ have the same functional properties as
$K^\pm(z)$, hence the morphism is well-defined.

The relations between $A_{R,\Lambda}(z),A_{R,\Lambda}(w)$
are such that
the relations between Cartan currents are preserved. One checks that
$$
\hbar[A_R(z),a(w)]
= \log q(z,w), \quad
\hbar[A_\Lambda(z),a(w)]
= \log q(w,q^{-\partial}z),
$$
$$
\hbar[A_R(z),b(w)]
= \log q(z,w),
\quad
\hbar[A_\Lambda(z),b(w)]
= \log q(w,z),
$$
(recall that $\on{log} q(z,w) = \sum_i \hbar(T+U)(\lambda_i)(z) r^i(w)$),
therefore the relations between the Cartan currents and the currents
$e(w)$, $f(w)$ are preserved.

These relations imply the relations
\begin{equation} \label{a:a}
[a(z),a(w)] = \on{log} q(z,w) + (1+q^{-\partial_z})
\on{log}{{w-z}\over{z-w}},
\end{equation}
\begin{equation} \label{a:b}
[b(z),b(w)] = -\on{log} q(z,w) + (1+q^{\partial_z})
\on{log}{{w-z}\over{z-w}},
\end{equation}
\begin{equation} \label{b:b}
[a(z),b(w)] = (1+q^{\partial_z})\on{log}{{w-z}\over{z-w}}.
\end{equation}

Here $\on{log}(w-z)/(z-w) = \on{log}(w) - \on{log}(z)
-\sum_{n\neq 0}(z/w)^n/n$.

Let us prove (\ref{a:a}).
We have
\begin{align*}
& \partial_z[a(z),a(w)] = {{\hbar\partial_z}\over{1-q^{-\partial_z}}}
[A_R(z) + A_\Lambda(z),a(w)]
= {{\partial_z}\over{1-q^{-\partial_z}}} \big( \on{log}q(z,w) +
\on{log}q(w,q^{-\partial}z)\big)
\\ & = {{\partial_z}\over{1-q^{-\partial_z}}} \big( \on{log}q(z,w) -
\on{log}q(q^{-\partial}z,w)\big)
+ {{\partial_z}\over{1-q^{-\partial_z}}} \big( \on{log} q(q^{-\partial}z,w)
+ \on{log} q(w,q^{-\partial}z)\big).
\end{align*}

Now
$$
\on{log} q(q^{-\partial} z,w)
+ \on{log} q(w,q^{-\partial}z) = \on{log} {{q^{-\partial}z-q^\partial w}
\over{z-w}} + \on{log}{{w-z}\over{q^{\partial}w-q^{-\partial}z}}
= (1-q^{-2\partial_z})\on{log}{{w-z}\over{z-w}},
$$
so
$$
\partial_z[a(z),a(w)] = \partial_z \sk{ \on{log} q(z,w)
+ (1+q^{-\partial_z}) \on{log}{{w-z}\over{z-w}}},
$$
therefore the difference of both sides of (\ref{a:a}) has the form $f(z)$.
Since this difference is also antisymmetric in $z,w$, it is zero.

Then (\ref{a:b}) follows from the addition of (\ref{a:a}) with
the formula expressing $[\hbar A_+(z),b(w)]$. The derivation of
(\ref{b:b}) is similar.

Let us now show that relation (\ref{e-e}) is preserved.
We have
$$
E(z)E(w) = \on{exp}(\alpha(z,w)) :E(z)E(w):,
$$
where $:E(z)E(w):\ := v(z)v(w)\on{exp}(a_+(z) + a_+(w))
\on{exp}(a_-(z) + a_-(w))$ and
$\alpha(z,w) := [a_+(z),a_-(w)] + {1\over 2}[a_+(z),a_+(w)]
+ {1\over 2}[a_-(z),a_-(w)]$.

Recall that we defined
$j(z,w) := {{z-q^\partial w}\over{q^{-\partial}z-w}}$,
$i(z,w) := q(z,w) (q^\partial z - w)/(z - q^\partial w)$.
Then
$$
[a(z),a(w)] = \on{log}(i(z,w)) +
\log{{z-q^\pa w}\over{q^\pa z - w}} + \log{{w-z}\over{z-w}}
+ \on{log}{{w-q^{-\pa}z}\over{q^{-\pa}z-w}}.
$$
Therefore
\begin{align*}
[a_+(z),a_+(w)] &
= \log(i)_{++}(z,w) - (\log(q^\pa z - w))_{++}
- (\log(q^{-\pa}z-w))_{++}
\\ &
= \log(i)_{++}(z,w) - \sk{\log{{q^\pa z - w}\over{z-q^{-\pa}w}}}_{++}
- \sk{\log{{q^{-\pa}z-w}\over{q^\pa z-w}}}_{++}
\\ & = \log(i)_{++}(z,w) + \log(j/j^{21})_{++}(z,w),
\end{align*}
\begin{align*}
[a_-(z),a_-(w)] & = \log(i)_{--}(z,w)
+ (\log(z-q^\pa w))_{--} + (\log(w-q^{-\pa}z))_{--}
\\ &
= \log(i)_{--}(z,w)
+ \sk{\log{{z-q^\pa w}\over{q^{-\pa}z-w}}}_{--}
+ \sk{\log{{w-q^{-\pa}z}\over{q^\pa w - z}}}_{--}
\\ & = \log(i)_{--}(z,w),
\end{align*}
\begin{align*}
[a_+(z),a_-(w)] & = \log(i)_{+-}(z,w)
+ \log(w-z) + (\log(w-z^{-\pa}z))_{+-}
\\ &
= \log(i)_{+-}(z,w) - (\log(w-q^{-\pa}z))_{--}
+ \log(w-z) + \log(w-z^{-\pa}z)
\\ &
= \log(i)_{+-}(z,w) - \sk{\log{{w-q^{-\pa}z}\over{q^\pa w - z}}}_{--}
+ \log(w-z) + \log(w-z^{-\pa}z)
\\ & = \log(i)_{+-}(z,w) + \log(j)_{--}(z,w)
+ \log(w-z)(w-q^{-\partial}z).
\end{align*}
Here we decompose an element $\beta$ of (a completion of)
$\widetilde\cK^{\otimes 2}$ as $\sum_{\eps,\eps'=+,-}\beta_{\eps,\eps'}$
according to the decomposition of $\cK$.

Finally,
\begin{align*}
\alpha(z,w) & = \Big( \log(i)_{+-} + {1\over 2}\log(i)_{++} + {1\over 2}
\log(i)_{--} + \log(j)_{--} + {1\over 2} \log(j/j^{21}) \Big)(z,w)
+ \log(w-z)(w-q^{-\pa}z)
\\ & = \gamma(z,w) + \log(w-z)(w-q^{-\pa}z),
\end{align*}
where $\gamma(z,w)\in \hbar \CC[[z,w]][z^{-1},w^{-1}][[\hbar]]$.
One checks that $\on{exp}(\gamma(z,w) - \gamma(w,z))
= i(z,w)j(z,w)/j(w,z)$: this follows from the fact that
$\log(j/j^{21})_{+-} = 0$, because
$$
\log(j/j^{21})(z,w) = \log(z-q^\pa w)(z-q^{-\pa}w)
- \log(q^{-\pa}z-w)(q^\pa z-w).
$$
Now we have
\begin{align*}
& (w-q^\partial z)E(z) E(w) = e^{\gamma(z,w)} (w-q^\partial z)(w-z)
(w-q^{-\partial}z):E(z)E(w):
\\ & = i(z,w)(q^{\partial}w-z)
e^{\gamma(w,z)} (w-z)
(q^{-\partial}w-z):E(z)E(w):
= i(z,w)(q^{\partial}w-z)
 E(w)E(z)
\end{align*}
(using the fact that $:E(z)E(w): = :E(w)E(z):$),
so relation (\ref{e-e}) is preserved.

One proves similarly that relation (\ref{f-f}) is preserved.

Let us now show that (\ref{e-f}) is preserved.
We have $E(z) F(w) = \on{exp}(\lambda(z,w)):E(z) F(w):$,
where $:E(z)F(w): = v(z)v'(w)\on{exp}(a_-(z)-b_-(w))\on{exp}(a_+(z)-b_+(w))$
and $\lambda(z,w) = -[a_+(z),b_-(w)]-{1\over 2}[a_+(z),b_+(w)]
- {1\over 2} [a_-(z),b_-(w)]$, and
$F(w)E(z) = \on{exp}(\mu(z,w)):E(z)F(w):$, where
$\mu(z,w) = -[b_+(w),a_-(z)]-{1\over 2}[b_+(w),a_+(z)]
-{1\over 2}[b_-(w),a_-(z)]$. So
$$
[E(z),F(w)] = (\on{exp}(\lambda(z,w)) - \on{exp}(\mu(z,w))):E(z)F(w):.
$$
Now
$$
[a(z),b(w)] = \log(w-z) + \log(w - q^\pa z) - \log(z-w) - \log(q^\pa z - w),
$$
Calculating as above, we find
$$
[a_+(z),b_+(w)] = -\log(j^{21})_{++}(z,w),
$$
$$
[a_+(z),b_-(w)] = \log(w-z)(w-q^\pa z) - \log(j^{21})_{--}(z,w)
$$
$$
[a_-(z),b_+(w)] = -\log(z-w)(q^\pa z-w) + \log(j^{21})_{++}(z,w),
$$
$$
[a_-(z),b_-(w)] = \log(j^{21})_{--}(z,w).
$$
Therefore
$$
\lambda(z,w) =
\Big( {1\over 2}\log(j^{21})_{++} + {1\over 2}\log(j^{21})_{--}\Big)(z,w)
- \log(w-z)(w-q^\pa z)
$$
and
$$
\mu(z,w) =
\Big( {1\over 2}\log(j^{21})_{++} + {1\over 2}\log(j^{21})_{--}\Big)(z,w)
- \log(z-w)(q^\pa z-w).
$$
Recall that
$\phi(z,w) = \Big( {1\over 2}\log(j^{21})_{++}
+ {1\over 2}\log(j^{21})_{--}\Big)(z,w)$. Then
\begin{align*}
& [E(z),F(w)] = e^{\phi(z,w)} \sk{ {1\over {(w-z)(w-q^\partial z)}}
- {1\over {(z-w)(q^{\partial}z-w)}}
} :E(z)F(w):
\\ & = {{e^{\phi(z,w)}}\over{z-q^{\partial} z}} \big( \delta(w-z)
- \delta(w-q^{\partial} z)\big)  :E(z)F(w):
\\ & =
{{r_0(w)}\over{(z-q^{\partial}z)/\hbar}}
{1\over\hbar}\sk{ e^{\phi(z,z)}\delta(z,w)
:E(z)F(z): - e^{\phi(q^{-\pa}w,w)}\delta(q^\pa z,w)
:E(q^{-\partial} w)F(w):} .
\end{align*}
Here $\delta(z-w) = 1/(z-w)+1/(w-z)$. We have $\delta(z,w) =
\delta(z-w)/r_0(w)$ (recall that $\omega(z) = r_0(z) dz$).

Now
\begin{align*}
& :E(z)F(z): = v(z)v'(z) e^{(a-b)_-(z)}e^{(a-b)_+(z)}
= v(z)v'(z)e^{\hbar (A_R)_-(z)}
e^{\hbar (A_R)_+(z)}
\\ & = v(z)v'(z)
e^{{1\over 2}[\hbar (A_R)_-(z),\hbar (A_R)_+(z)]}
u_R(z)^{-1}K^+(z).
\end{align*}
Now $[A_R(z),A_R(w)] = 0$, hence $[(A_R)_+(z),(A_R)_-(w)] = 0$,
which implies $:E(z)F(z): = v(z)v'(z) K^+(z)$.

On the other hand,
\begin{align*}
& :E(q^{-\partial}w)F(w): = v(q^{-\pa w})v'(w) e^{(q^{-\partial}a-b)_-(w)}
e^{(q^{-\partial}a-b)_+(w)} = v(q^{-\pa w})v'(w) e^{-\hbar (A_\Lambda)_-(w)}
e^{-\hbar (A_\Lambda)_+(w)}
\\ & = v(q^{-\pa w})v'(w) e^{{1\over 2}[\hbar (A_\Lambda)_-(w),\hbar
(A_\Lambda)_+(w)]}u_\Lambda(w)K_-(w)^{-1}.
\end{align*}

Now $[\hbar A_\Lambda(z),\hbar A_\Lambda(w)] = (q^{-\partial}\otimes
q^{-\partial}-1)\on{log}q(z,w)\in \hbar^2 \wedge^2(R)[[\hbar]]$.
Then we get
$$
:E(q^{-\partial}z)F(z): = v(q^{-\pa}z)v'(z) e^{( (q^{-\partial}\otimes
q^{-\partial}-1)\on{log}q(z,w))_{-+|w=z}} u_\Lambda(z)K_-(z)^{-1},
$$
with
$\big( (q^{-\partial}\otimes
q^{-\partial}-1)\on{log}q(z,w)\big)_{-+|w=z}\in \hbar^2\cK[[\hbar]]$.

Finally,
\begin{align*}
& [E(z),F(w)] =
v(z) v'(w){{r_0(w)}\over{(z-q^{\partial}z)/\hbar}}\times
\\ & {1\over \hbar}\sk{ \delta(z,w)
e^{\phi(z,z)} u_R(z)^{-1}K^+(z) - \delta(z,q^{-\partial} w)
e^{\phi(q^{-\pa}w,w)}e^{( (q^{-\partial}\otimes
q^{-\partial}-1)\on{log}q(w,z))_{-+|z=w}} u_\Lambda(w)K_-(w)^{-1}
} .
\end{align*}

Recall that $\alpha(z) = \phi(z,z)$ and $\beta(w) = \phi(q^{-\pa}w,w)
+ \Big((q^{-\pa}\otimes q^{-\pa}-1)\log q(w,z)\Big)_{-+|z=w}$.

Then
\begin{align*}
& \skk{{{(z-q^\pa z)/\hbar}\over{v(z)}} E(z),
e^{-\alpha_\Lambda(w) - \beta_R(w)}{{r_0(w)}\over{v'(w)}}F(w)} =
\\ &\quad = {1\over \hbar}\sk{ \delta(z,w)
e^{\alpha_R(z)} u_R(z)^{-1}K^+(z)
- \delta(z,q^{-\partial} w)
e^{\beta_\Lambda(w)}
u_\Lambda(w)K_-(w)^{-1}
} .
\end{align*}

Now since $v(z) = (z-q^\pa z)/\hbar$, $v'(w) =
r_0(w)e^{-\alpha_\Lambda(w)-\beta_R(w)}$, $u_R(z) = e^{\alpha_R(z)}$
and $u_\Lambda(w) = e^{-\beta_\Lambda(w)}$, relation (\ref{e-f}) is preserved.
\hfill \qed \medskip

\section{Double of elliptic quantum group
$E_{\tau,\eta}(\mathfrak{sl}_2)$ versus
elliptic algebra $U_{q,p}(\widehat{\mathfrak{sl}}_2)$ }

It was shown in \cite{EF98} that the algebra given by the commutation
relations \r{e-e} to \r{e-f} with the function
$q(z,w)=\frac{\theta(z-w+\hbar)}{\theta(z-w-\hbar)}$
given in Subsubsection~\ref{2.4.3}
is isomorphic to the centrally extended double of the elliptic quantum group
$E_{\tau,\eta}(\mathfrak{sl}_2)$.
This isomorphism is an analogue of the Ding-Frenkel isomorphism \cite{DF}.
The basic representations considered here provide level one representations of this
central extended double of the elliptic quantum group.

Let us write a basis in $\mathcal{K}$ in the case of the data given in
Subsubsection~\ref{2.4.3}. For $\la\neq 0$ it is
\bigskip

\noindent
\begin{tabular}{c|cccccccc}
\hline $i$ &  $\cdots$ & -3 & -2 & -1 & 0 & 1 & 2 & $\cdots$ \\ \hline
&&&&&&&&\\
$\epsilon^{i;\la}(z)$ &  $\cdots$ & $z^2$ & $z^1$ & $z^0$ &
$\frac{\theta(z+\la)}{\theta(z)\theta(\la)}$ &
$\frac{1}{1!}\left(\frac{\theta(z+\la)}{\theta(z)\theta(\la)}\right)^\prime$ &
$\frac{1}{2!}\left(\frac{\theta(z+\la)}{\theta(z)\theta(\la)}\right)^{\prime\prime}$ & $\cdots$
\\ 
&&&&&&&&\\
$\epsilon_{i;\la}(w)$ &  $\cdots$ &
$\frac{1}{2!}\left(\frac{\theta(\la+w)}{\theta(w)\theta(\la)}\right)^{\prime\prime}$ &
$\frac{-1}{1!}\left(\frac{\theta(\la+w)}{\theta(w)\theta(\la)}\right)^{\prime}$ &
$\frac{\theta(\la+w)}{\theta(w)\theta(\la)}$ & $(-w)^0$ & $(-w)^1$ & $(-w)^2$ & $\cdots$ \\
&&&&&&&&\\ \hline
\end{tabular}
\bigskip

\noindent
and for $\la=0$
\bigskip

\noindent
\begin{tabular}{c|cccccccc}
\hline
$i$ &  $\cdots$ & -3 & -2 & -1 & 0 & 1 & 2 & $\cdots$ \\ \hline
&&&&&&&&\\
$\epsilon^{i;0}(z)$ &  $\cdots$ & $z^2$ & $z^1$ & $z^0$ &
$(\ln\theta(z))'$ &
$\frac{1}{1!}\left(\ln\theta(z)\right)^{\prime\prime}$ &
$\frac{1}{2!}\left(\ln\theta(z)\right)^{\prime\prime\prime}$ & $\cdots$
\\ 
&&&&&&&&\\
$\epsilon_{i;0}(w)$ &  $\cdots$ &
$\frac{1}{2!}\left(\ln\theta(w)\right)^{\prime\prime\prime}$ &
$\frac{-1}{1!}\left(\ln\theta(w)\right)^{\prime\prime}$ &
$(\ln\theta(w))'$ & $(-w)^0$ & $(-w)^1$ & $(-w)^2$ & $\cdots$ \\
&&&&&&&&\\ \hline
\end{tabular}
\bigskip

The total currents $e(z)$ and $f(z)$ can be decomposed into modes with respect to
the basis given in the first table, namely
\begin{equation}\label{to-ef-de}
e(z)=\sum_{i\in\ZZ} e[i;\la]\epsilon^{i;\lambda}(z)\,,\qquad
f(z)=\sum_{i\in\ZZ} f[i;\la]\epsilon^{i;\lambda}(z)
\end{equation}
Using a duality of the bases in $\mathcal{K}$:
$\oint dw\ \epsilon^{i;\la}(w)\ \epsilon_{j;\la}(w)=\delta_{i,j}$,
we can write the modes as the contour integrals
\begin{equation}\label{ef-mod}
e[i,\la]=\oint_{C^+_0} dw\ e(w)\ \epsilon_{i;\la}(w)\,,\qquad
f[i,\la]=\oint_{C^+_0} dw\ f(w)\ \epsilon_{i;\la}(w)
\end{equation}
where $C^+_0$ is a contour encircling 0 counterclockwise. Define half currents as the sums
\begin{equation}\label{h-cur}
e^\mp_\la(z)=\mp\sum_{i\geq 0\atop i<0} e[i;\mp\la] \epsilon^{i;\la}(z)\,,\qquad
f^\mp_\mu(z)=\mp\sum_{i\geq 0\atop i<0} f[i;\mp\mu] \epsilon^{i;\mu}(z)
\end{equation}
we may easily deduce that they can be written in the integral forms
\begin{equation}\label{ha-e}
e^-_\la(z)=\oint_{C_0^-}\frac{\theta(u-z-\la)}{\theta(u-z)\theta(-\la)}\ e(u)\ du\,,\qquad
e^+_\la(z)=\oint_{C_{0,z}^+}\frac{\theta(u-z-\la)}{\theta(u-z)\theta(-\la)}\ e(u)\ du
\end{equation}
and
\begin{equation}\label{ha-f}
f^-_{-\la}(z)=\oint_{C_0^-}\frac{\theta(u-z+\la)}{\theta(u-z)\theta(\la)}\ f(u)\ du\,,\qquad
f^+_{-\la}(z)=\oint_{C_{0,z}^+}\frac{\theta(u-z+\la)}{\theta(u-z)\theta(\la)}\ f(u)\ du
\end{equation}
with $C^+_{0,z}$ a contour encircling the points $0$ and $z$
counterclockwise and
$C^-_0$ a contour encircling $0$ clockwise. In order to obtain
positive and negative half-currents \r{ha-e} and \r{ha-f} we used
the Taylor expansion of the kernel
$\frac{\theta(z-w+\la)}{\theta(z-w)\theta(\la)}$ in powers of $z$
and $w$ respectively.
Note that the expansion in powers of $z$ is defined in the domain
$z\ll w$ and in the domain
$w\ll z$ with respect to powers of $w$. This explains the
difference between the contours of integrations
in \r{ha-e} and \r{ha-f}. The same reason is behind the appearance
of the $\delta$-function
in the sum
\begin{equation}\label{del-dec}
\sum_{i\geq 0} \epsilon^{i,\la}(z)\ \epsilon_{i,\la}(w)+\sum_{i<0}
\epsilon^{i,-\la}(z)\ \epsilon_{i,-\la}(w)=
\frac{\theta(z-w+\la)}{\theta(z-w)\theta(\la)}+
\frac{\theta(w-z-\la)}{\theta(w-z)\theta(\la)}=\delta(z-w)
\end{equation}

Besides the Taylor expansion, the kernel $\frac{\theta(z-w+\la)}{\theta(z-w)\theta(\la)}$
has the following Fourier mode
expansion\footnote{An analogous formula exists for the elliptic current
algebra associated with Baxter's 8-vertex $R$-matrix, see \cite{KLPST}.}
\begin{equation}\label{Fur}
\frac{\theta(z+\la)}{\theta(z)\theta(\la)}
=\sum_{n\in\ZZ}\frac{e^{2\pi \i z(\la+n)/\tau }}{e^{2\pi \i (\la+n)/\tau }-1}=
\sum_{n\in\ZZ}\frac{e^{2\pi \i \la(z+n)/\tau }}{e^{2\pi \i (z+n)/\tau }-1}
\end{equation}
which is valid for $0<\Im(z/\tau)<\Im(1/\tau)$ and $0<\Im(\la/\tau)<\Im(1/\tau)$.
Formula \r{Fur} leads to the temptation to write instead of \r{h-cur}
half-currents in the form
\begin{equation}\label{h-cur-F}
e^\pm_\la(z)=\sum_{n\in\ZZ}\ e[n] \
\frac{e^{2\pi \i z(n+\la)/\tau}}{1-e^{\pm 2\pi \i(n+\la)/\tau}}\,,\qquad
f^\pm_\mu(z)=\sum_{n\in\ZZ}\ f[n] \
\frac{e^{2\pi \i z(n+\mu)/\tau}}{1-e^{\pm 2\pi \i(n+\mu)/\tau}}
\end{equation}
where $e[n]$ are $f[n]$ are $\la$-independent modes of the total currents
\begin{equation*}
e(z)= \sum_{n\in\ZZ}\ e[n] \ e^{2\pi \i z n/\tau}\,,\qquad
f(z)= \sum_{n\in\ZZ}\ f[n] \ e^{2\pi \i z n/\tau}
\end{equation*}
such that
\begin{equation*}
e(z)= (e^+_{\la}(z)+e^-_{\la}(z)) \ e^{-2\pi \i z \la/\tau}\,,\qquad
f(z)= (f^+_{\mu}(z)+f^-_{\mu}(z)) \ e^{-2\pi \i z \mu/\tau}
\end{equation*}

The definition of the currents of type \r{h-cur-F}
was used in the paper \cite{JKOS} to build the
elliptic algebra $U_{q,p}(\widehat{\mathfrak{sl}}_2)$.
The resulting algebra of half-currents appears to be completely different from the algebra
investigated in \cite{EF98}. In particular, it acquires a dynamical behaviour not only with respect
to the parameter $\la$, but also with respect to the elliptic parameter $\tau$.
Moreover, because of the formal relation between positive and negative half-currents
of the form $e^+_\la(u+1)=-e^-_\la(u)$ only one $L$-operator, say $L_\la^+(u)$, is sufficient
to describe the whole algebra. All these phenomena were described in detail
in the paper \cite{KLP-CMP} in the example of the scaling limit of the elliptic algebra.
A classical version of this scaled elliptic algebra was investigated in the paper \cite{KLPST} where,
in particular, it was proved that on the classical level, the central extension of the classical
currents should be done by the cocycle which is defined
by a differentiation with respect to moduli of the elliptic curve.
The decomposition into Fourier modes becomes in the scaling limit
the decomposition into continuous Fourier modes and the summations on
$n\in\ZZ$ become integrals along the real
axis. The dynamical behaviour with respect to the parameter $\la$
disappears in the scaling limit.
The resulting algebra is a current algebra on the infinite cylinder
and it was used in the papers \cite{KLP-CMP,CKP} to explain
the dynamical symmetries of the Sine-Gordon model.

This situation, when different decompositions of the total currents
lead to different
algebras for half-currents, was investigated also in the paper \cite{KLP-AMS}
in the case of the centrally extended Yangian double. Elliptic algebras
considered in the papers \cite{EF98} and \cite{JKOS} are in fact elliptic
generalizations of two constructions of the centrally
extended Yangian double given in \cite{KLP-AMS}.

\section{Conclusion}

To conclude we would like to mention that the existence of the level zero representations
in both elliptic algebras formulated in the papers \cite{EF98} and \cite{JKOS}
is based on the Fay  identity
\begin{equation*}
\begin{split}\theta(a+c)\;\theta(b+d)\theta(a-c)\;\theta(b-d)&=
\theta(a+b)\;\theta(c+d)\theta(a-b)\;\theta(c-d) \ny\\ &-
\theta(a+d)\;\theta(c+b)\theta(a-d)\;\theta(c-b)
\end{split}
\end{equation*}
which can be written in the form \cite{JKOS}
\begin{equation*}
\frac{\theta(u_1+t)}{\theta(u_1)}\ \eta_{s,t}(u_2)=
 \frac{\theta(u_1-u_2+t)}{\theta(u_1-u_2)}\ \eta_{s+t,t}(u_2)+
 \eta_{s,t}(u_2-u_1)\eta_{s+t,t}(u_1)
\end{equation*}
where
\begin{equation*}
\eta_{s,t}(u)= \frac{\theta(u+s)\theta(t)}{\theta(u)\theta(s)}
\end{equation*}

The Fay identity can be written for higher genus algebraic curve and it is interesting
to apply this identity for investigation of higher genus current algebras. It is clear that
formulas like \r{Fur} does not exist in higher genus, so the decompositions of total currents
which lead to dynamical shifts in modules of the curve are impossible. But ``non-symmetric''
decompositions of the type \r{h-cur} are always possible.
It is interesting to combine the Fay identity and these decompositions to investigate some current
algebras associated with higher genus algebraic curves.
The work in this direction is in progress.

\section*{Acknowledgements}
The research by S.P. was supported in part by the grants INTAS-OPEN-03-51-3350,
Heisenberg-Landau program, RFBR grant 03-02-17373.
V.R. greatly acknowledged the support of the  RFBR grant 03-02-17554
and of franco-ukrainian Collaboration Program ``Dnipro'' during
his visit to Kiev. He thanks Bogolyubov ITP of NAN Ukraine for
hospitality. S.P. and V.R. are grateful to the RFBR grant
to support scientific schools  NSh-1999.2003.2.

This work was
started when S.P. enjoyed a CNRS visitor position
in Angers. He acknowledges the CNRS support and the Mathematical
Departments of the Angers and Strasbourg Universities for
hospitality. Part of the work was done during
the visit of two authors (S.P. and V.R.) at the Max Planck Institut
f\"ur Mathematik (Bonn). They wish to thank the Institut
for hospitality and for its stimulating scientific atmosphere.

B.E. would also like to thank M. Jimbo for explanations about the
relations between \cite{EF98} and \cite{JKOS}.

\end{document}